\documentclass[12pt]{article}
\usepackage{amscd,amsfonts,amsmath,amssymb,yfonts}
\usepackage{amstext,amsthm,mathbbol,mathrsfs,stmaryrd} 
\usepackage[cmtip,arrow]{xy}
\usepackage{pb-diagram,pb-xy} 
\usepackage{curves,graphicx,texdraw} 
\usepackage{tikz} 
\usetikzlibrary{calc} 
\usepackage[pagewise,switch]{lineno} 

\theoremstyle{definition} 
\newtheorem{definition}{\sc Definition}[section] 
\theoremstyle{plain}
\newtheorem{theorem}{\sc Theorem}[section]

\newtheorem{lemma}[theorem]{\sc Lemma} 
\newtheorem{corollary}[theorem]{\sc Corollary} 
\newcounter{conj} 
\newcounter{prob} 
\newcounter{ques} 
\newcounter{exmpl} 
\newcounter{exer} 
\def\exmpl{\vskip 1.1\baselineskip\noindent
\textsc{Example\ }\theexmpl : 
\stepcounter{exmpl}} 
 
\begin{document}

\title{$k$-uniform $l$-Eulerian hypergraphs}

\author{
Yiming Lei\thanks{Corresponding author. Email address: leiyiming008@gmail.com (Yiming Lei).},
Chenxi Lv, Bing Dai\\
\small College of Mathematical Sciences, Bohai University, Jinzhou, 121000, China
}

\date{}
\maketitle

\begin{abstract}
We consider $l$-Euler trails and $l$-Euler circuits in finite $k$-uniform directed multihypergraphs, where $k\geq 2$ and $1\leq l<k$. Let $D$ be such a directed multihypergraph. The notions of $l$-connectivity and $l$-balance are defined in terms of the $l$-tails and $l$-heads of its hyperarcs. It is shown that $D$ is $l$-Eulerian if and only if $D$ is $l$-connected and $l$-balanced. Moreover, $D$ is $l$-semi-Eulerian if and only if $D$ is $l$-connected and has exactly two states with unequal $l$-indegrees and $l$-outdegrees, one with $l$-outdegree one greater than its $l$-indegree and the other with $l$-indegree one greater than its $l$-outdegree. The proofs are based on an $l$-transition digraph and the classical Euler theorem for directed multigraphs. The corresponding conditions for $k$-uniform multihypergraphs are also obtained in terms of their orientations. For a given orientation, an expected linear-time algorithm is provided to determine whether an $l$-Euler trail or an $l$-Euler circuit exists and to construct one when it does.
\end{abstract}

\vskip 2mm
\noindent\textbf{Keywords}:
Directed multihypergraph, Eulerian hypergraph, Transition digraph, Overlap trail

\vskip 2mm
\noindent\textbf{Mathematics Subject Classification}:
05C45, 05C65

\section{Hypergraph}

Let $k\geq 2$ and $1\leq l<k$.
In this paper, all hypergraphs and directed hypergraphs are finite.
A \textit{$k$-uniform multihypergraph} $H=(V,E)$ consists of a nonempty
vertex set $V$ and a finite multiset $E$ of $k$-element subsets of $V$.
Each hyperedge contains $k$ distinct vertices, but the same $k$-element
subset may appear more than once in $E$.
If two hyperedges have the same vertex set, we regard them as different
hyperedge occurrences in $E$. 

Let $P=v_1,v_2,\ldots,v_m$ be a sequence of vertices with 
$m=k+q(k-l)$ for some integer $q\geq 0$. 
If for every $0\leq i\leq q$ there is a hyperedge occurrence $e_i\in E(H)$ 
whose underlying $k$-set is 
$\{v_{i(k-l)+1},v_{i(k-l)+2},\ldots,v_{i(k-l)+k}\}$, 
and the occurrences $e_0,e_1,\ldots,e_q$ are pairwise distinct, 
then $P$ is called an \textit{$l$-trail} of $H$. 
A cyclic sequence $C=(v_1,v_2,\ldots,v_m)$ with 
$m=q(k-l)$ for some positive integer $q$ is called 
an \textit{$l$-circuit} of $H$ if, with subscripts read modulo $m$, 
for every $0\leq i\leq q-1$ there is a hyperedge occurrence $e_i\in E(H)$ 
whose underlying $k$-set is 
$\{v_{i(k-l)+1},v_{i(k-l)+2},\ldots,v_{i(k-l)+k}\}$, 
and the occurrences $e_0,e_1,\ldots,e_{q-1}$ are pairwise distinct. 
Let $E(P)$ and $E(C)$ denote the multisets of hyperedge occurrences appearing in 
$P$ and $C$, respectively. 
An \textit{$l$-Euler trail} of $H$ is an $l$-trail $P$ such that 
$E(P)=E(H)$. 
An \textit{$l$-Euler cycle} of $H$ is an $l$-circuit $C$ such that 
$E(C)=E(H)$. 
The hypergraph $H$ is called \textit{$l$-Eulerian} if it has an $l$-Euler cycle, 
and \textit{$l$-semi-Eulerian} if it has an $l$-Euler trail but no $l$-Euler cycle. 

Note that in the definitions of $l$-trails and $l$-circuits,
the overlap of two consecutive hyperedges is determined by $l$
consecutive positions in the vertex sequence.
We do not require that the intersection of their vertex sets
has exactly $l$ elements.
This distinction will be used in the construction of the
$l$-transition digraph. 

In \cite{Kuhn:2010}, Hamilton $l$-cycles in uniform hypergraphs
were studied.
Several notions of Euler tours in hypergraphs have also been
considered; see
\cite{Lonc:2010,Lonc:2017,Glock:2020,SajnaWagner:2023,LoPigaSanhueza:2024}.
The construction of the transition digraph in this paper is related
to the auxiliary digraphs used for universal cycles and overlap cycles;
see \cite{ChungDiaconisGraham:1992,Jackson:1993,HoranHurlbert:2013}. 

A \textit{$k$-uniform directed multihypergraph} $D=(V,E)$ consists
of a nonempty vertex set $V$ and a finite multiset $E$ of hyperarcs.
Each hyperarc occurrence $a\in E$ is represented by an ordered
$k$-tuple $a=(v_1,v_2,\ldots,v_k)$ of distinct vertices.
Different hyperarc occurrences may have the same ordered $k$-tuple. 
For a hyperarc $a=(v_1,v_2,\ldots,v_k)$, the sequence 
$t_l(a)=(v_1,v_2,\ldots,v_l)$ is called the \textit{$l$-tail} of $a$, 
and the sequence $h_l(a)=(v_{k-l+1},\allowbreak v_{k-l+2},\ldots,v_k)$ is called 
the \textit{$l$-head} of $a$. 
Let $S$ be an ordered $l$-tuple of vertices of $D$.
The number of hyperarc occurrences whose $l$-tail is $S$ is denoted by
$d_D^{l+}(S)$, and the number of hyperarc occurrences whose $l$-head
is $S$ is denoted by $d_D^{l-}(S)$. That is,
\[
 d_D^{l+}(S)=\bigl|\{a\in E(D):t_l(a)=S\}\bigr|,\qquad
 d_D^{l-}(S)=\bigl|\{a\in E(D):h_l(a)=S\}\bigr|.
\]
Here we count hyperarcs with their multiplicities. 

A sequence of hyperarc occurrences $P=(a_1,a_2,\ldots,a_r)$ is called 
an \textit{$l$-trail} of $D$ if no hyperarc occurrence is repeated and 
$h_l(a_i)=t_l(a_{i+1})$ for every $1\leq i<r$. 
An $l$-trail $C=(a_1,a_2,\ldots,a_r)$ is called 
an \textit{$l$-circuit} of $D$ if $h_l(a_r)=t_l(a_1)$. 
An \textit{$l$-Euler trail} of $D$ is an $l$-trail containing 
every hyperarc occurrence of $D$ exactly once. 
An \textit{$l$-Euler circuit} of $D$ is an $l$-circuit containing 
every hyperarc occurrence of $D$ exactly once. 
The directed hypergraph $D$ is called \textit{$l$-Eulerian} 
if it has an $l$-Euler circuit, and \textit{$l$-semi-Eulerian} 
if it has an $l$-Euler trail but no $l$-Euler circuit. 

Let $D$ be a directed multihypergraph.
By replacing every hyperarc occurrence of $D$ with its underlying
vertex set, we obtain a multihypergraph $H(D)$, called the
\textit{underlying multihypergraph} of $D$.
The multiplicities of the hyperarcs are preserved.
Conversely, let $H$ be a $k$-uniform multihypergraph.
If we choose an ordering of the $k$ vertices of each hyperedge
occurrence of $H$, then the resulting directed multihypergraph $D$
is called an \textit{orientation} of $H$. 

\begin{lemma}
Let $H$ be a $k$-uniform multihypergraph. 
Then $H$ has an $l$-Euler trail if and only if 
there exists an orientation $D$ of $H$ that has an $l$-Euler trail. 
Moreover, $H$ has an $l$-Euler cycle if and only if 
there exists an orientation $D$ of $H$ that has an $l$-Euler circuit. 
\label{lemma:orientation}
\end{lemma}

\begin{proof}
Suppose that $H$ has an $l$-Euler trail. 
Let $P=v_1,v_2,\ldots,v_m$ be an $l$-Euler trail of $H$, where 
$m=k+q(k-l)$ for some integer $q\geq 0$. 
For each $0\leq i\leq q$, let $e_i$ be the hyperedge occurrence of $P$ whose 
underlying $k$-set is 
$\{v_{i(k-l)+1},v_{i(k-l)+2},\ldots,v_{i(k-l)+k}\}$. 
Since $P$ is an $l$-Euler trail, $e_0,e_1,\ldots,e_q$
are pairwise distinct and $E(P)=E(H)$. 
For each $0\leq i\leq q$, orient $e_i$ as the hyperarc occurrence 
$a_i=(v_{i(k-l)+1},v_{i(k-l)+2},\ldots,v_{i(k-l)+k})$. 
Hence these orderings define an orientation $D$ of $H$. 
For every $0\leq i<q$, the last $l$ entries of $a_i$ are
$v_{(i+1)(k-l)+1},\ldots,v_{(i+1)(k-l)+l}$.
These are also the first $l$ entries of $a_{i+1}$. 
Hence $h_l(a_i)=t_l(a_{i+1})$ for every $0\leq i<q$. 
Thus $(a_0,a_1,\ldots,a_q)$ is an $l$-trail of $D$. 
Since it contains every hyperarc occurrence of $D$ exactly once, 
it is an $l$-Euler trail of $D$. 

Next, suppose that $H$ has an $l$-Euler cycle. 
Let $C=(v_1,v_2,\ldots,v_m)$ be an $l$-Euler cycle of $H$, where 
$m=q(k-l)$ for some positive integer $q$, and read the subscripts modulo $m$. 
Orient the hyperedge occurrence corresponding to the $i$th block as 
$a_i=(v_{i(k-l)+1},v_{i(k-l)+2},\ldots,v_{i(k-l)+k})$ for 
$0\leq i\leq q-1$. 
By the definition of these hyperarcs, we have
$h_l(a_i)=t_l(a_{i+1})$, where the subscripts of the hyperarcs
are read modulo $q$. 
Hence $(a_0,a_1,\ldots,a_{q-1})$ is an $l$-circuit containing every hyperarc 
occurrence exactly once. 
Therefore, the resulting orientation $D$ has an $l$-Euler circuit. 

Suppose, conversely, that an orientation $D$ of $H$ has an $l$-Euler trail 
$P=(a_1,a_2,\ldots,a_r)$. 
Write $a_i=(x_{i,1},x_{i,2},\ldots,x_{i,k})$ for each $1\leq i\leq r$. 
Since $P$ is an $l$-trail, we have $h_l(a_i)=t_l(a_{i+1})$ for every 
$1\leq i<r$. 
Thus $x_{i,k-l+j}=x_{i+1,j}$ for every $1\leq i<r$ and $1\leq j\leq l$. 
Start with the sequence $x_{1,1},x_{1,2},\ldots,x_{1,k}$. 
For $i=2,3,\ldots,r$, append the last $k-l$ entries 
$x_{i,l+1},x_{i,l+2},\ldots,x_{i,k}$ of $a_i$. 
Since $x_{i,k-l+j}=x_{i+1,j}$ for $1\leq j\leq l$,
the $k$ consecutive entries beginning at position $(i-1)(k-l)+1$
are $x_{i,1},x_{i,2},\ldots,x_{i,k}$ for every $1\leq i\leq r$.
Hence the corresponding hyperedge occurrence of $H$ is the
underlying hyperedge occurrence of $a_i$. 
Since $P$ contains every hyperarc occurrence of $D$ exactly once and $D$ is an 
orientation of $H$, the resulting vertex sequence contains every hyperedge occurrence 
of $H$ exactly once. 
Hence it is an $l$-Euler trail of $H$. 

Finally, suppose that $P=(a_1,a_2,\ldots,a_r)$ is an $l$-Euler circuit of $D$. 
We construct a vertex sequence as above.
Since $h_l(a_r)=t_l(a_1)$, the last $l$ entries of this sequence
are the same as its first $l$ entries.
We delete the last $l$ entries and consider the remaining
$r(k-l)$ entries as a cyclic sequence. 
Then the $k$ consecutive cyclic entries corresponding to each $a_i$ form the 
underlying hyperedge occurrence of $a_i$, including the last block across the cyclic 
boundary. 
All hyperedge occurrences of $H$ appear exactly once. 
Therefore, this cyclic sequence is an $l$-Euler cycle of $H$. 
\end{proof}

\section{Eulerian Hypergraph}

\begin{definition}
Let $D$ be a $k$-uniform directed multihypergraph. 
The \textit{$l$-state set} of $D$ is 
$\mathcal{S}_l(D)=\{t_l(a):a\in E(D)\}\cup\{h_l(a):a\in E(D)\}$. 

We say that $D$ is \textit{$l$-connected} if for any two states 
$S,T\in\mathcal{S}_l(D)$, there exist states 
$S_0,S_1,\ldots,S_r$ with $S_0=S$ and $S_r=T$ such that every two 
consecutive states are the $l$-tail and the $l$-head of a hyperarc occurrence 
of $D$.
We say that $D$ is \textit{$l$-balanced} if 
$d_D^{l-}(S)=d_D^{l+}(S)$ for every $S\in\mathcal{S}_l(D)$. 
\end{definition}

We now consider the conditions for a directed multihypergraph
to be $l$-Eulerian or $l$-semi-Eulerian.
We first define a directed multigraph which will be used in the proofs. 

\begin{definition}
Let $D$ be a $k$-uniform directed multihypergraph. 
The \textit{$l$-transition digraph} of $D$, denoted by 
$\mathcal{T}_l(D)$, is the directed multigraph with vertex set 
$V(\mathcal{T}_l(D))=\mathcal{S}_l(D)$ such that each hyperarc occurrence 
$a\in E(D)$ corresponds to one directed edge 
$e_a:t_l(a)\longrightarrow h_l(a)$. 
By this definition, every hyperarc occurrence of $D$
corresponds to exactly one edge of $\mathcal{T}_l(D)$, and conversely.
Hence the correspondence
$E(D)\longleftrightarrow E(\mathcal{T}_l(D))$
is a bijection preserving multiplicities. 

A directed multigraph is called \textit{weakly connected} if its underlying 
undirected multigraph is connected.
\end{definition}

\begin{lemma}
The directed hypergraph $D$ is $l$-connected if and only if 
$\mathcal{T}_l(D)$ is weakly connected. 
Moreover, for every $S\in\mathcal{S}_l(D)$, 
$d_{\mathcal{T}_l(D)}^+(S)=d_D^{l+}(S)$ and 
$d_{\mathcal{T}_l(D)}^-(S)=d_D^{l-}(S)$. 
\label{lemma:transitionproperties}
\end{lemma}

\begin{proof}
Suppose that $D$ is $l$-connected. 
Let $S,T\in V(\mathcal{T}_l(D))=\mathcal{S}_l(D)$. 
By the definition of $l$-connectivity, there exist states 
$S=S_0,S_1,\ldots,S_r=T$ such that for every $1\leq i\leq r$ there is a 
hyperarc occurrence $a_i\in E(D)$ with 
$\{t_l(a_i),h_l(a_i)\}=\{S_{i-1},S_i\}$. 
By the definition of $\mathcal{T}_l(D)$, the edge $e_{a_i}$
joins $S_{i-1}$ and $S_i$ in the underlying undirected graph.
Hence $S_0,S_1,\ldots,S_r$ is a walk connecting $S$ and $T$.
Since $S,T\in\mathcal{S}_l(D)$ were arbitrary,
$\mathcal{T}_l(D)$ is weakly connected. 

Suppose, conversely, that $\mathcal{T}_l(D)$ is weakly connected. 
Let $S,T\in\mathcal{S}_l(D)$. 
There is a walk $S=S_0,S_1,\ldots,S_r=T$ in the underlying graph of 
$\mathcal{T}_l(D)$. 
For every $1\leq i\leq r$, the vertices $S_{i-1}$ and $S_i$ are the endpoints of 
some edge $e_{a_i}$ corresponding to a hyperarc occurrence $a_i\in E(D)$. 
Therefore, $\{t_l(a_i),h_l(a_i)\}=\{S_{i-1},S_i\}$ for every $1\leq i\leq r$. 
By the definition of $l$-connectivity, $D$ is $l$-connected. 

Now let $S\in\mathcal{S}_l(D)$. 
By the definition of $\mathcal{T}_l(D)$, an edge $e_a$ leaves $S$
if and only if $t_l(a)=S$.
Since every hyperarc occurrence corresponds to exactly one edge,
we have $d_{\mathcal{T}_l(D)}^+(S)=d_D^{l+}(S)$.
Similarly, an edge $e_a$ enters $S$ if and only if $h_l(a)=S$.
Hence $d_{\mathcal{T}_l(D)}^-(S)=d_D^{l-}(S)$. 
\end{proof}

\begin{lemma}
Let $a_1,a_2,\ldots,a_r$ be hyperarc occurrences of $D$, and let $Q$ be the 
sequence of their corresponding edges in $\mathcal{T}_l(D)$, taken in the same order. 
Then $(a_1,a_2,\ldots,a_r)$ is an $l$-trail in $D$ if and only if $Q$ is 
a directed trail in $\mathcal{T}_l(D)$. 
Moreover, $(a_1,a_2,\ldots,a_r)$ is an $l$-circuit in $D$ if and only if $Q$ 
is a closed directed trail in $\mathcal{T}_l(D)$. 
\label{lemma:trailcorrespondence}
\end{lemma}

\begin{proof}
Suppose that $(a_1,a_2,\ldots,a_r)$ is an $l$-trail in $D$. 
Since no hyperarc occurrence is repeated, the corresponding edges 
$e_{a_1},e_{a_2},\ldots,e_{a_r}$ are pairwise distinct. 
For every $1\leq i<r$, we have $h_l(a_i)=t_l(a_{i+1})$. 
By the definition of the transition digraph, $e_{a_i}$ ends
at $h_l(a_i)$ and $e_{a_{i+1}}$ starts at $t_l(a_{i+1})$.
Since $h_l(a_i)=t_l(a_{i+1})$, these two edges are consecutive
in a directed trail for every $1\leq i<r$. 
Hence $Q=(e_{a_1},e_{a_2},\ldots,e_{a_r})$ is a directed trail in 
$\mathcal{T}_l(D)$. 

Suppose, conversely, that $Q=(e_{a_1},e_{a_2},\ldots,e_{a_r})$ is a directed 
trail in $\mathcal{T}_l(D)$. 
The edges of $Q$ are pairwise distinct. 
By the multiplicity-preserving bijection between $E(D)$ and 
$E(\mathcal{T}_l(D))$, the hyperarc occurrences 
$a_1,a_2,\ldots,a_r$ are also pairwise distinct. 
Since $Q$ is a directed trail, the terminal vertex of $e_{a_i}$
is the initial vertex of $e_{a_{i+1}}$ for every $1\leq i<r$. 
Hence $h_l(a_i)=t_l(a_{i+1})$ for every $1\leq i<r$. 
Thus $(a_1,a_2,\ldots,a_r)$ is an $l$-trail in $D$. 

If $(a_1,a_2,\ldots,a_r)$ is an $l$-circuit, then in addition 
$h_l(a_r)=t_l(a_1)$. 
Therefore, the terminal vertex of $e_{a_r}$ is the initial vertex of $e_{a_1}$, 
and $Q$ is a closed directed trail. 
Suppose, conversely, that $Q$ is closed.
Then the last edge of $Q$ ends at the initial vertex of its first edge.
Hence $h_l(a_r)=t_l(a_1)$. 
Hence $(a_1,a_2,\ldots,a_r)$ is an $l$-circuit. 
\end{proof}

\begin{corollary}
The directed hypergraph $D$ has an $l$-Euler trail 
if and only if $\mathcal{T}_l(D)$ has an Euler trail. 
Moreover, $D$ has an $l$-Euler circuit if and only if 
$\mathcal{T}_l(D)$ has an Euler circuit. 
\label{corollary:eulercorrespondence}
\end{corollary}

We will use the classical Euler theorem for directed multigraphs.
A finite directed multigraph with at least one edge has an Euler circuit
if and only if all its non-isolated vertices are connected in the
underlying undirected graph and each vertex has equal indegree and outdegree.
It has an Euler trail but no Euler circuit if and only if all its
non-isolated vertices are connected in the underlying undirected graph,
there is exactly one vertex whose outdegree is one larger than its indegree,
there is exactly one vertex whose indegree is one larger than its outdegree,
and all other vertices have equal indegree and outdegree. 

\begin{theorem}
Let $D$ be a $k$-uniform directed multihypergraph with $E(D)\neq\emptyset$. 
Then $D$ is $l$-Eulerian if and only if $D$ is $l$-connected and $l$-balanced. 
\label{theorem:directedeuler}
\end{theorem}

\begin{proof}
Suppose that $D$ is $l$-Eulerian. 
Let $C=(a_1,a_2,\ldots,a_r)$ be an $l$-Euler circuit of $D$.
Then $C$ contains every hyperarc occurrence of $D$ exactly once. 
By Lemma \ref{lemma:trailcorrespondence}, 
$Q=(e_{a_1},e_{a_2},\ldots,e_{a_r})$ is a closed directed trail in 
$\mathcal{T}_l(D)$. 
Since $C$ contains all hyperarc occurrences of $D$ exactly once
and each hyperarc occurrence corresponds to one edge of $\mathcal{T}_l(D)$,
$Q$ contains every edge of $\mathcal{T}_l(D)$ exactly once. 
Hence $Q$ is an Euler circuit of $\mathcal{T}_l(D)$. 

By the definition of $\mathcal{S}_l(D)$, every state is the
$l$-tail or the $l$-head of some hyperarc occurrence of $D$.
Hence $\mathcal{T}_l(D)$ has no isolated vertices. 
Since $\mathcal{T}_l(D)$ has an Euler circuit, the classical Euler theorem implies 
that $\mathcal{T}_l(D)$ is weakly connected and that 
$d_{\mathcal{T}_l(D)}^-(S)=d_{\mathcal{T}_l(D)}^+(S)$ for every 
$S\in\mathcal{S}_l(D)$. 
By Lemma \ref{lemma:transitionproperties}, $D$ is $l$-connected and 
$d_D^{l-}(S)=d_D^{l+}(S)$ for every $S\in\mathcal{S}_l(D)$. 
Thus $D$ is $l$-connected and $l$-balanced. 

Suppose, conversely, that $D$ is $l$-connected and $l$-balanced. 
By Lemma \ref{lemma:transitionproperties}, $\mathcal{T}_l(D)$ is weakly connected. 
Moreover, for every $S\in\mathcal{S}_l(D)$, we have 
$d_{\mathcal{T}_l(D)}^-(S)=d_D^{l-}(S)=d_D^{l+}(S)=d_{\mathcal{T}_l(D)}^+(S)$. 
Since $E(D)\neq\emptyset$, $\mathcal{T}_l(D)$ has at least one edge.
By the definition of $\mathcal{S}_l(D)$, it has no isolated vertices. 
Therefore, the classical Euler theorem implies that $\mathcal{T}_l(D)$ has an 
Euler circuit. 
Let $Q=(e_{a_1},e_{a_2},\ldots,e_{a_r})$ be such an Euler circuit. 
Since $Q$ is an Euler circuit, every edge of $\mathcal{T}_l(D)$
appears in $Q$ exactly once.
Hence the corresponding hyperarc occurrences $a_1,a_2,\ldots,a_r$
are all the hyperarc occurrences of $D$, each appearing exactly once. 
By Lemma \ref{lemma:trailcorrespondence}, 
$(a_1,a_2,\ldots,a_r)$ is an $l$-circuit in $D$. 
Therefore, it is an $l$-Euler circuit of $D$ and $D$ is $l$-Eulerian. 
\end{proof}

\begin{theorem}
Let $D$ be a $k$-uniform directed multihypergraph with $E(D)\neq\emptyset$. 
Then $D$ is $l$-semi-Eulerian if and only if 
$D$ is $l$-connected and there exist two distinct states 
$S^+,S^-\in\mathcal{S}_l(D)$ such that
\[
\begin{aligned}
d_D^{l+}(S^+)&=d_D^{l-}(S^+)+1,\\
d_D^{l-}(S^-)&=d_D^{l+}(S^-)+1,
\end{aligned}
\]
and $d_D^{l-}(S)=d_D^{l+}(S)$ for every 
$S\in\mathcal{S}_l(D)\setminus\{S^+,S^-\}$. 
In this case, the $l$-tail of the first hyperarc of every $l$-Euler trail 
is $S^+$, and the $l$-head of the last hyperarc is $S^-$. 
\label{theorem:directedsemieuler}
\end{theorem}

\begin{proof}
Suppose that $D$ is $l$-semi-Eulerian. 
Then $D$ has an $l$-Euler trail but has no $l$-Euler circuit. 
Let $P=(a_1,a_2,\ldots,a_r)$ be an $l$-Euler trail of $D$. 
By Lemma \ref{lemma:trailcorrespondence}, 
$Q=(e_{a_1},e_{a_2},\ldots,e_{a_r})$ is a directed trail in 
$\mathcal{T}_l(D)$. 
Since $P$ contains every hyperarc occurrence of $D$ exactly once, $Q$ contains 
every edge of $\mathcal{T}_l(D)$ exactly once. 
Thus $Q$ is an Euler trail of $\mathcal{T}_l(D)$. 
If $\mathcal{T}_l(D)$ had an Euler circuit, Corollary 
\ref{corollary:eulercorrespondence} would imply that $D$ had an $l$-Euler circuit, 
which is impossible. 
Therefore, $\mathcal{T}_l(D)$ has an Euler trail but no Euler circuit. 

Let $S^+=t_l(a_1)$ and $S^-=h_l(a_r)$. 
By the definition of $\mathcal{T}_l(D)$, $S^+$ is the initial
vertex of $Q$ and $S^-$ is its terminal vertex.
Since $D$ has no $l$-Euler circuit, $Q$ is not closed.
Hence $S^+\neq S^-$. 
By the classical Euler theorem for directed multigraphs, 
$\mathcal{T}_l(D)$ is weakly connected, 
$d_{\mathcal{T}_l(D)}^+(S^+)=d_{\mathcal{T}_l(D)}^-(S^+)+1$, 
$d_{\mathcal{T}_l(D)}^-(S^-)=d_{\mathcal{T}_l(D)}^+(S^-)+1$, and every other 
state has equal indegree and outdegree. 
By Lemma \ref{lemma:transitionproperties}, $D$ is $l$-connected
and $S^+$ and $S^-$ satisfy the required degree conditions. 

Suppose, conversely, that $D$ is $l$-connected and there exist two distinct 
states $S^+,S^-\in\mathcal{S}_l(D)$ satisfying the stated degree conditions. 
By Lemma \ref{lemma:transitionproperties}, $\mathcal{T}_l(D)$ is weakly connected. 
The same lemma gives 
$d_{\mathcal{T}_l(D)}^+(S^+)=d_{\mathcal{T}_l(D)}^-(S^+)+1$ and 
$d_{\mathcal{T}_l(D)}^-(S^-)=d_{\mathcal{T}_l(D)}^+(S^-)+1$, while every other 
vertex of $\mathcal{T}_l(D)$ has equal indegree and outdegree. 
Therefore, by the classical Euler theorem, $\mathcal{T}_l(D)$ has an Euler trail 
from $S^+$ to $S^-$. 
Since $S^+\neq S^-$, this Euler trail is not closed.
Moreover, since the indegree and outdegree of $S^+$ are different,
$\mathcal{T}_l(D)$ has no Euler circuit. 
By Corollary \ref{corollary:eulercorrespondence}, $D$ has an $l$-Euler trail and 
has no $l$-Euler circuit. 
Hence $D$ is $l$-semi-Eulerian. 

We now determine the first and last $l$-states of an $l$-Euler trail.
In any Euler trail of $\mathcal{T}_l(D)$, the initial vertex
has outdegree one larger than indegree, and the terminal vertex has
indegree one larger than outdegree.
Since $S^+$ and $S^-$ are the unique vertices with these properties,
every Euler trail of $\mathcal{T}_l(D)$ starts at $S^+$ and ends at $S^-$. 
By Lemma \ref{lemma:trailcorrespondence}, the $l$-tail of the first hyperarc of 
every $l$-Euler trail of $D$ is $S^+$, and the $l$-head of its last hyperarc is $S^-$. 
\end{proof}

\begin{theorem}
Let $H$ be a $k$-uniform multihypergraph with $E(H)\neq\emptyset$. 
Then $H$ is $l$-Eulerian if and only if 
there exists an orientation $D$ of $H$ that is $l$-connected and $l$-balanced. 
\label{theorem:undirectedeuler}
\end{theorem}

\begin{proof}
Suppose that $H$ is $l$-Eulerian. 
Then $H$ has an $l$-Euler cycle. 
By Lemma \ref{lemma:orientation}, there exists an orientation $D$ of $H$ that 
has an $l$-Euler circuit. 
Thus $D$ is $l$-Eulerian. 
By Theorem \ref{theorem:directedeuler}, $D$ is $l$-connected and $l$-balanced. 
Hence $H$ has an orientation that is $l$-connected and $l$-balanced. 

Suppose, conversely, that $H$ has an orientation $D$
which is $l$-connected and $l$-balanced. 
By Theorem \ref{theorem:directedeuler}, $D$ is $l$-Eulerian.
Hence $D$ has an $l$-Euler circuit. 
Since $D$ is an orientation of $H$, Lemma \ref{lemma:orientation} implies that 
$H$ has an $l$-Euler cycle. 
Therefore, $H$ is $l$-Eulerian. 
\end{proof}

\begin{theorem}
Let $H$ be a $k$-uniform multihypergraph with $E(H)\neq\emptyset$. 
Then $H$ is $l$-semi-Eulerian if and only if $H$ has an $l$-connected 
orientation $D$ with exactly two distinct states $S^+,S^-\in\mathcal{S}_l(D)$ 
such that
\[
d_D^{l+}(S^+)=d_D^{l-}(S^+)+1,\qquad
d_D^{l-}(S^-)=d_D^{l+}(S^-)+1,
\]
and $d_D^{l-}(S)=d_D^{l+}(S)$ for every 
$S\in\mathcal{S}_l(D)\setminus\{S^+,S^-\}$, while $H$ has no 
$l$-connected and $l$-balanced orientation.
\label{theorem:undirectedsemieuler}
\end{theorem}

\begin{proof}
Suppose that $H$ is $l$-semi-Eulerian. 
Then $H$ has an $l$-Euler trail but has no $l$-Euler cycle. 
By Lemma \ref{lemma:orientation}, there exists an orientation $D$ of $H$ having 
an $l$-Euler trail. 
Suppose that $D$ has an $l$-Euler circuit.
By Lemma \ref{lemma:orientation}, $H$ has an $l$-Euler cycle.
This contradicts the assumption that $H$ is $l$-semi-Eulerian.
Therefore, $D$ has no $l$-Euler circuit, and hence $D$ is
$l$-semi-Eulerian. 
By Theorem \ref{theorem:directedsemieuler}, $D$ is $l$-connected and there exist 
two distinct states $S^+,S^-\in\mathcal{S}_l(D)$ such that 
$d_D^{l+}(S^+)=d_D^{l-}(S^+)+1$, 
$d_D^{l-}(S^-)=d_D^{l+}(S^-)+1$, and 
$d_D^{l-}(S)=d_D^{l+}(S)$ for every 
$S\in\mathcal{S}_l(D)\setminus\{S^+,S^-\}$. 

We now show that $H$ has no $l$-connected and $l$-balanced
orientation.
Suppose that $H$ has such an orientation $D'$.
By Theorem \ref{theorem:directedeuler}, $D'$ has an $l$-Euler circuit.
By Lemma \ref{lemma:orientation}, $H$ has an $l$-Euler cycle.
This contradicts the assumption that $H$ is $l$-semi-Eulerian.
Hence $H$ has no $l$-connected and $l$-balanced orientation. 

Suppose, conversely, that $H$ has an $l$-connected orientation $D$ with two 
distinct states $S^+,S^-\in\mathcal{S}_l(D)$ such that 
$d_D^{l+}(S^+)=d_D^{l-}(S^+)+1$, 
$d_D^{l-}(S^-)=d_D^{l+}(S^-)+1$, and 
$d_D^{l-}(S)=d_D^{l+}(S)$ for every 
$S\in\mathcal{S}_l(D)\setminus\{S^+,S^-\}$, and suppose that $H$ has no 
$l$-connected and $l$-balanced orientation. 
By Theorem \ref{theorem:directedsemieuler}, $D$ has an $l$-Euler trail. 
By Lemma \ref{lemma:orientation}, $H$ also has an $l$-Euler trail. 

We now show that $H$ has no $l$-Euler cycle.
Suppose that $H$ has an $l$-Euler cycle. 
By Lemma \ref{lemma:orientation}, there exists an orientation $D'$ of $H$ that 
is $l$-Eulerian. 
By Theorem \ref{theorem:directedeuler}, $D'$ is $l$-connected and 
$l$-balanced. 
This contradicts the assumption that $H$ has no $l$-connected and 
$l$-balanced orientation. 
Hence $H$ has no $l$-Euler cycle. 
Therefore, $H$ is $l$-semi-Eulerian. 
\end{proof}

\begin{corollary}
Let $k$ and $l$ be fixed and suppose that an orientation $D$
is given.
Under the standard RAM model with expected constant-time
dictionary operations, we can determine whether $D$ is $l$-Eulerian
or $l$-semi-Eulerian in expected time $O(|E(D)|)$.
If an $l$-Euler trail or an $l$-Euler circuit exists,
we can also construct it in expected time $O(|E(D)|)$. 
\label{corollary:algorithm}
\end{corollary}

\begin{proof}
Let $m=|E(D)|$.
We first scan all hyperarc occurrences of $D$ once.
For each hyperarc occurrence $a=(v_1,v_2,\ldots,v_k)$,
we compute its $l$-tail $t_l(a)$ and its $l$-head $h_l(a)$.
Since $k$ and $l$ are fixed, this takes constant time.
We use a dictionary to assign an index to every state that appears
as an $l$-tail or an $l$-head, and insert the directed edge
$t_l(a)\longrightarrow h_l(a)$ into $\mathcal{T}_l(D)$.
We also count the outdegree of $t_l(a)$ and the indegree of $h_l(a)$.
Since every hyperarc has one $l$-tail and one $l$-head,
there are at most $2m$ states.
Hence this construction takes expected time $O(m)$ under the stated
assumption on dictionary operations. 

We then consider the underlying undirected graph of
$\mathcal{T}_l(D)$ and perform a graph search.
By Lemma \ref{lemma:transitionproperties}, this determines whether
$D$ is $l$-connected in time $O(m)$.
From the indegrees and outdegrees, we can determine in the same time
whether every state is balanced, or whether exactly two states satisfy
the degree conditions in Theorem \ref{theorem:directedsemieuler}.
Hence, by Theorems \ref{theorem:directedeuler} and
\ref{theorem:directedsemieuler}, the recognition algorithm
takes expected time $O(m)$. 

Suppose that the conditions for an $l$-Euler circuit hold.
We apply Hierholzer's algorithm to $\mathcal{T}_l(D)$
and obtain an Euler circuit.
If the conditions for an $l$-Euler trail but no $l$-Euler circuit
hold, we use the open-trail version of the same algorithm,
starting at the state whose outdegree is one larger than its indegree.
In either case, each edge is processed a constant number of times.
Hence the construction takes time $O(m)$.
We then replace every edge $e_a$ of the resulting Euler trail or circuit
by its corresponding hyperarc occurrence $a$.
By Lemma \ref{lemma:trailcorrespondence}, we obtain an $l$-Euler trail
or an $l$-Euler circuit of $D$.
Therefore, the total expected running time is $O(|E(D)|)$. 
\end{proof}

\exmpl
We now give an example to show that the $l$-connectivity of $D$
cannot be replaced by a connectivity condition on the underlying
hypergraph $H(D)$.
Let $k=3$ and $l=1$, and suppose that $D$ has the six hyperarcs 
$(1,7,2)$, $(2,8,3)$, $(3,9,1)$, $(4,7,5)$, $(5,10,6)$, and $(6,11,4)$. 
Every state has equal $1$-indegree and $1$-outdegree. 
The hyperedges of $H(D)$ are connected by intersections
of cardinality one: the first three hyperedges form one $1$-cycle,
the last three form another, and the hyperedges $\{1,7,2\}$
and $\{4,7,5\}$ have the common vertex $7$. 
However, the six $1$-states split into the two components 
$\{1,2,3\}$ and $\{4,5,6\}$ under the $1$-connectivity relation. 
Hence $D$ is not $1$-connected and therefore is not $1$-Eulerian. 
Hence this example shows that connectivity of the underlying
hypergraph does not imply $l$-connectivity of $D$. 

\section{Concluding remarks}

By Theorem \ref{theorem:directedeuler}, a directed multihypergraph
is $l$-Eulerian if and only if it is $l$-connected and $l$-balanced.
This is similar to the classical Euler theorem for directed multigraphs.
In this paper, we use the $l$-transition digraph to reduce the problem
of finding $l$-Euler trails and $l$-Euler circuits in a given orientation
to the problem of finding Euler trails and Euler circuits in a
directed multigraph.
For a multihypergraph $H$ without a given orientation, we need to
determine whether there exists an $l$-connected and $l$-balanced
orientation of $H$.
We may further consider conditions expressed directly in terms of
invariants of $H$, and the computational complexity of finding
such an orientation. 

We may also consider the case in which two consecutive hyperedges
are required to have exactly $l$ common vertices.
In this case, the equality of the $l$-head of one hyperarc and
the $l$-tail of the next hyperarc does not exclude other common
vertices of these hyperarcs.
Hence the transition digraph alone is not sufficient for this problem.
Further conditions on the transitions would be required in addition
to the degree conditions considered above. 

\section*{Acknowledgements}\label{section: acknowledgements} 
This work is supported by Youth project of basic scientific research program of Liaoning provincial department of education: LJ212510167015, 
Ph.D. research startup foundation of 
Bohai University: 0522bs006 and 
Research foundation of 
Bohai University: 0522xn031.

\end{document}